\documentclass[12pt]{amsart}

\usepackage{tikz}
\usepackage{xcolor}
\usepackage{skull}
\usepackage{comment}
\usepackage{enumitem}
\usepackage{fullpage,multirow, amsmath, amsthm,amsfonts,amssymb,stmaryrd, mathrsfs,hhline}
\usepackage{mathdots}
\usepackage{mathrsfs}
\usepackage{mathtools}

\newtheorem{theorem}{Theorem}[section]
\newtheorem{lemma}[theorem]{Lemma}
\newtheorem{prop}[theorem]{Proposition}

\theoremstyle{definition}
\newtheorem{definition}[theorem]{Definition}
\newtheorem{example}[theorem]{Example}
\newtheorem{remark}[theorem]{Remark}
\newtheorem{notation}[theorem]{Notation}
\newtheorem{corollary}[theorem]{Corollary}
\newtheorem*{theorem*}{Theorem}

\numberwithin{equation}{section}

\newcommand{\tha}{^{\text{th}}}
\newcommand{\OO}{\mathcal O}
\newcommand{\pp}{\mathfrak{p}}   

\newcommand{\PP}{\mathbb{P}}

\DeclareMathOperator{\PGL}{PGL}
\DeclareMathOperator{\wdeg}{wdeg}
\DeclareMathOperator{\good}{good}
\DeclareMathOperator{\an}{an}

\usepackage{multirow}

\begin{document}

\title{Local fields, iterated extensions, and Julia Sets}

\author{Pui Hang Lee, Michelle Manes, Nha Xuan Truong}

\address{Pui Hang Lee, Mathematics department, University of Hawaii at Manoa, 2565 McCarthy Mall (Keller Hall 401A)
Honolulu, Hawaii 96822, USA}
\email{plee9@hawaii.edu}

\address{Michelle Manes, American Institute of Mathematics,
 Caltech 8-32
 1200 E California Blvd
 Pasadena CA 91125, USA}
 \email{mmanes@aimath.org}

\address{Nha Xuan Truong, Beijing International Center for Mathematical Research, Peking University, 5 Yi He Yuan Road, Haidian District, Beijing, 100871, China. }
\email{nxtruong@bicmr.pku.edu.cn}

\begin{abstract}
Let $K$ be a field complete with respect to a discrete valuation $v$ of residue characteristic $p$. For $\alpha \in K$, let $K_\infty$ be the extension obtained by adjoining all iterated preimages of $\alpha$ under a unicritical polynomial $f_c(z)=z^\ell - c \in K[z]$. We study the extension $K_\infty/K$ and show that its qualitative behavior depends only on the valuation of $c$. This removes the previous restrictions on $\ell$ in work of Anderson--Hamblen--Poonen--Walton and completes the classification for all $\ell \ge 2$. We also relate the ramification to the structure of the Berkovich Julia set of $f_c$.
\end{abstract}

\maketitle

\section{Introduction}
Let $K$ be a local field and let $\overline{K}$ be a fixed separable closure. Let $g(z)\in K[z]$ be a separable polynomial of degree $\ell \geq 2$. We define the  $n\tha$ iterate of $g$ by the $n$-fold composition: $g^n(z) = g\circ g \circ \cdots \circ g (z)$. 
Fix $\alpha \in K$ and define the set $V_n = \{ \beta \in \overline{K} : g^n(\beta) = \alpha \}$. We then define a tower of algebraic extensions of $K$ by 
\[
K_n := K\left(V_n\right)
\quad \text{and} \quad
K_\infty := \bigcup_{n\geq0}K_n.
\]

The current work is motivated by ~\cite{anderson-hamblen-poonen-walton}, in which the authors studied the fields arising from this construction for unicritical polynomials $f_c(z) = z^\ell -c$ where either $p\nmid \ell$ or $p = \ell$.
 The authors identify a cutoff value 
\[
 \nu_\infty = 
     -\frac{\ell}{\ell-1} v(\ell) 
\]
such that:
\begin{itemize}
    \item for $v(c) < \nu_\infty$, the extension $K_\infty / K$ is finite,
    \item for $v(c) = \nu_\infty$, the extension is infinite but can be finitely ramified, depending on the valuation of the root point $v(\alpha)$, and
    \item for $v(c) > \nu_\infty$, the extension is infinite with infinite ramification (infinite wild ramification if $p=\ell$).
    
\end{itemize} 
Our initial goal was to extend these results to all possible degrees $\ell \geq 2$. In the cases considered in~\cite{anderson-hamblen-poonen-walton}, we observed that $v(c) \geq \nu_\infty$ if and only if the polynomial $f_c(z) = z^\ell -c$ has potential good reduction. Moreover, this valuation cutoff corresponds precisely to a change in the structure of the Berkovich Julia set: when $v(c) < \nu_\infty$, the Berkovich Julia set for $f_c$ is a Cantor set of Type I points, while for $v(c) \geq \nu_\infty$, it is a single type II point. We suspected this dichotomy might be related to the extension $K_\infty / K$, and that the Berkovich Julia set could explain this relationship. However, when we consider all degrees $\ell \geq 2$, the situation turns out to be more complicated.

Writing $\ell = Np^k$ with $(p,N) = 1$, we define a second important cutoff value:
\[
 \nu_{\text{good}} = 
 \begin{cases}
     0 & \text{if } N > 1 \\
     -\frac{p}{p-1} 
     & \text{otherwise.}
 \end{cases}
\]
We are now ready to state our main theorem.
\begin{theorem*}
Suppose that $\ell \geq 2$, $(\ell,p) \neq 1$, and $c \in \overline{K}$. 
Consider the polynomial $f_c(z) = z^\ell - c$. Then
\begin{enumerate}
    \item $f_c$ has potential good reduction if and only if $v(c) \geq \nu_{\text{good}}$.
    \item (Berkovich Julia set of $f_c$):
    \begin{itemize}
        \item If $v(c) < \nu_\infty$, the Berkovich Julia set is a Cantor set 
        consisting entirely of type I points.
        \item If $\nu_\infty \leq v(c) < \nu_{\text{good}}$, the Berkovich Julia set is a Cantor set containing 
        type II points and possibly type IV points.
        \item If $v(c) \geq \nu_{\text{good}}$, the Berkovich Julia set consists of a single type II point.
    \end{itemize}
    \item (Ramification of $K_\infty/K$):
    \begin{itemize}
        \item If $v(c) < \nu_\infty$, then $K_\infty/K$ is a finite extension.
        \item If $v(c) = \nu_\infty$, then $K_\infty/K$ is an infinite extension. 
        It is finitely ramified if and only if $\ell = p$ and $\alpha$ lies within 
        the closed unit disk centered at a fixed point of $f$.
        \item If $v(c) > \nu_\infty$, then $K_\infty/K$ is infinitely wildly ramified.
    \end{itemize}
    
    \end{enumerate}
\end{theorem*}
\begin{remark}\begin{enumerate}
    \item The theorem recovers the result in~\cite{anderson-hamblen-poonen-walton} regarding the absence of deep ramification, which differs from the expectation in~\cite{aitken-hajir-maire} that preimage trees of a generic polynomial of degree divisible by $p$ should exhibit deep ramification.
    
    \item Note that $\nu_\infty \leq \nu_{\text{good}}$, with equality when $p \nmid \ell$ or $p =\ell$ (exactly the cases considered in~\cite{anderson-hamblen-poonen-walton}). This explains the clean transition we observed before.

    \item In the boundary case $v(c) = \nu_\infty$, the method in \cite{anderson-hamblen-poonen-walton} applies only for $\ell = p$, whereas our approach works for all cases except $\ell = p$, providing a unified result for all $\ell$.
\end{enumerate}
\end{remark}
We hope that this more dynamical viewpoint may allow future work on other interesting families beyond the unicritical case. Especially in families where the Julia set mixes type I, II, or IV points, more intricate arboreal Galois behavior and more interesting dynamics should arise, raising largely open questions about the ramification of iterated extensions. A geometric approach to local arboreal problems may prove essential for addressing such questions, offering insights less accessible through purely algebraic methods.
\subsection*{Outline of the paper}

In Section~\ref{sec:potgoodred}, we define potential good reduction for polynomials and determine the necessary and sufficient condition for $f_c(z) = z^\ell - c$ to have potential good reduction. In Section~\ref{sec:NP}, we study the Newton polygon of the polynomial $f(z) = (z+y)^\ell - y^\ell - d$, which allows us to relate the valuations $v(x-y)$ and $v(f(x)-f(y))$ for points $x$ and $y$ in the backward orbit of $\alpha$ under $f_c$. Section~\ref{sec:JuliaSets} describes the Berkovich Julia sets of the polynomials $f_c(z)$. Finally, in Section~\ref{sec:extension}, we use these dynamical tools to prove our results on ramification.

\section*{Acknowledgements}
The authors thank Rob Benedetto for suggesting the direction of this research. The second author was supported by the Simons Foundation Travel Support for Mathematicians (grant 963056). This article is part of the Ph.D. thesis of the first author. The authors also express gratitude to the contributors of the SAGE software, as some of our arguments were initially tested through computer simulations.

\section{Criteria for potential good reduction}\label{sec:potgoodred}
We begin with some notation.  

 \begin{tabular}{ll}
 $K$ & a non-archimedean field, complete with respect to a valuation $v$.\\
$ \overline{K}$ & a fixed separable closure of $K$.\\
$\OO_K$ & the ring of integers $\{ \alpha \in K \colon v(\alpha) \geq 0\}$.\\
$\pp$ & the maximal ideal of $\OO_K$.\\
$\overline a$ & the image of $a \in \OO_K$ under the natural reduction map $\OO_K \to \OO_K / \pp\OO_K$.\\
$\overline g$ & the image of a polynomial $g \in \OO_K[z]$ obtained by reducing the coefficients.
\end{tabular}

 Let $p$ be the characteristic of the residue field $\OO_K / \pp\OO_K$, and normalize the valuation so that $v(p) = 1$. 
The reduction map above induces a reduction map 
\begin{align*}
    \PP^1(K) &\to \PP^1(\OO_K / \pp\OO_K)\\
    [x:y] &\mapsto [\bar x : \bar y],
\end{align*}
for $x,y \in \OO_K$ chosen so that $\min\{v(x),v(y)\} = 0$. This lets us extend the reduction map to $K$ in a natural way, and points in $\PP^1(K) \smallsetminus \OO_K$ have image the point at infinity.
 
 Throughout this paper, we 
write 
\[f_c(z) = z^\ell - c \quad
 \text{ with } c \in  \overline{K} \text{ and } \ell \geq 2.
 \]
 Further we factor the exponent
 \[
 \ell = Np^k, \quad  \text{ with } (N,p) = 1, \   N\geq 1, \text{ and } k\geq 0.
 \]
 Finally, we define the two key ``cutoff values'' for $v(c)$ that are the basis of our investigation:
\begin{align}
 \nu_{\good} &= 
 \begin{cases}
     0 & \text{if } N > 1 \\
     -\frac{p}{p-1} 
     & \text{otherwise.}
 \end{cases}\\
 \nu_\infty &= 
     -\frac{\ell}{\ell-1} v(\ell).
\end{align}
In arithmetic geometry, the notion of ``good reduction'' refers to a variety under a reduction map having the same properties (for example genus of a curve) as the original. This motivates the definition of good reduction for dynamical systems.
\begin{definition}
\label{def:bad.red}
Let $f(z) \in \mathcal O_K[z]$.
We say that $f$ has {\bf good reduction} if $\deg(\overline f) = \deg(f)$, and $f$ has {\bf potential good reduction} if there is some $\phi(z) \in \PGL_2(\overline{K})$, such that $f^\phi := \phi^{-1}\circ f \circ \phi$ has good reduction. If $f$ does not have good reduction, then it has {\bf bad reduction}. If $f$ does not have potential good reduction, we say it has {\bf persistent bad reduction}.
\end{definition}
Though the definition considers only polynomials with coefficients $\OO_K$, it can be extended to polynomials with coefficients in the field $K$ by considering them as rational functions with homogeneous presentations. (We refer the reader to~\cite[Section~4.3]{Ben01} for details.) In the case $f_c(z) = z^\ell -c$, we conclude that $f_c$ has good reduction precisely when $v(c)\geq 0$.

We now give necessary and sufficient conditions, based on $\ell$ and the valuation of the constant term $c$, for $f_c$ to have potential good reduction. 
The following lemma of Benedetto allows us to characterize when a polynomial has potential good reduction.

\begin{lemma}[\protect{\cite[Corollary 4.6]{Ben01}}]
\label{lem:pgr} 

Let $f\in K[z]$ be a polynomial, and let $g$ be a polynomial conjugate of $f$ such that $g$ is monic and $g(0)=0$. Then $f$ has potential good reduction if and only if $g$ has  good reduction.
\end{lemma}

In other words, $f$ has potential good reduction if and only if it is conjugate to a monic polynomial that fixes 0 and has all coefficients with non-negative valuation.
Proposition~\ref{prop:goodreduction} makes this precise for unicritical polynomials $f_c$, connecting the valuation of $c$ being sufficiently negative with the map  having persistent bad reduction. The strategy is to move one of the fixed points of $f$ to $0$, resulting in a polynomial conjugate to $f$ and meeting the hypotheses of Benedetto's lemma.

We begin with a useful formula for the $p$-adic valuation of binomial coefficients, which will help us to analyze the coefficients of this conjugate polynomial.

\begin{lemma}
\label{val_binomial}
Let $\ell=Np^k$, where $p\nmid N$. For all $n\leq p^k$ we have $$v\left(\binom{\ell}{n}\right)=k-v(n).$$
\end{lemma}
\begin{proof}
Using the  formula for binomial coefficients, we see that 
$$ \binom{\ell}{n}= \frac{\ell}{n} \binom{\ell-1}{n-1}$$
It suffices to show that $v\left(\binom{\ell-1}{n-1}\right)=0.$ By Kummer's Theorem, $v\left(\binom{\ell-1}{n-1}\right)$ is the number of carries when adding $n-1$ and $\ell - n$ written in base $p$. We observe that in base $p$, the last $k$ digits of $\ell-1= Np^k-1$ are all $p-1$, so there are no such carries when $n \leq p^k$.  \qedhere
\end{proof}

\begin{prop}
\label{prop:goodreduction}
Let $f_c(z) = z^\ell-c \in \overline{K}[z]$ be a separable polynomial. 
Then  $f_c$ has potential good reduction if and only if 
    $v(c)\geq \nu_{\good}$.
\end{prop}

\begin{proof}
If $v(c)\geq 0$, then $f_c(z)$ already has good reduction, so  assume that $v(c) < 0$. Let $b$ be a fixed point of $f_c(z)$, so that $b^\ell-c=b.$ We conjugate $f_c(z)$ by $\phi(z)=z+b$ to  $\displaystyle g(z):=\phi^{-1}\circ f\circ \phi$, which satisfies $g(0)=0$ and is monic. Then $f_c(z)$ has potential good reduction if and only if $g(z)$ has good reduction. 
We have 
$$\displaystyle g(z)=\sum\limits_{n=1}^{\ell-1}\binom{\ell}{n}b^{n}z^{\ell-n}. $$
Thus, $g(z)$ has good reduction if and only if $v\left(\binom{\ell}{n}b^n\right)\geq0$ for all $n = 1,\dots,\ell-1$.
Since $v(c)<0$,  we have $
0 > v(c)=v(b^\ell-b) \geq \min\{\ell v(b), v(b)\}.$
We conclude that $v(b)<0$ and $v(c) = \ell v(b)$. Therefore, $g(z)$ has good reduction if and only if  
\begin{equation}
v(c)\geq-v\left(\binom{\ell}{n}\right)\frac{\ell}{n}= -\frac{k-v(n)}{n}\ell \quad \textrm{ for all } 1 \leq n \leq \ell-1.
\label{eqn:v(c)bound}
\end{equation}
We observe that  
\[
\max_{1\leq n \leq \ell-1}\left\{-\frac{k-v(n)}{n}\right\}
=
\begin{cases}
     -\frac{p}{p-1}
    &\text{ if } \ell = p^k, \text{ for } k \geq 1, \\
    0 & \text{ if } \ell = Np^k \text{ for } N>1 \text{ and } (N,p)=1.
\end{cases}
\]
Thus the condition in~\eqref{eqn:v(c)bound} is equivalent to 
$v(c) \geq
\nu_{\good}$.
\end{proof}

\section{Newton polygon of $(z+y)^\ell-y^\ell-d$}\label{sec:NP}

The tools in this section allow us to extend the ``general lemmas'' in~\cite[Section~3]{anderson-hamblen-poonen-walton} to exponents $\ell$ where $p \mid \ell$ but $p \neq \ell$. The idea is the same as in the earlier paper: We want to understand $v(x-y)$ when we know that $f_c(x)-f_c(y) = d$ for some $ d \in \overline{K}$. To do this, we consider the Newton polygon of the polynomial $F(z) = (z+y)^\ell - y^\ell -d$ and apply our results to $z = x-y$.

The following technique will help our analysis of the Newton polygon. 
Suppose that we have a list of points: $$P_i = (i,y_i) \quad\textrm{with} \quad y_i \in \mathbb{R}, \textrm{ for } i = 0,1, \dots, n.$$ We may \emph{shift them relative to a linear function} $y = ax+b$ (for some $a,b \in \mathbb{R}$), by transforming them into $$Q_i  = (i, y_i + ai+b) \textrm{\quad for \quad}i = 0,1, \dots, n.$$
We observe that the $x$-coordinates of the break-points of the lower convex hulls of the $P_i$ and the $Q_i$ are exactly the same.  In particular: 
\begin{enumerate}
\item For a fixed $i_0$, the point $P_{i_0}$ is a vertex of the convex hull of the $P_i$ if and only if the point $Q_{i_0}$ is a vertex of the convex hull of the $Q_i$'s.
\item Let $m_{AB}$ represent the slope of the segment from point $A$ to point $B$. Then for every $i\neq j$ we have $m_{Q_iQ_j} = m_{P_iP_j}+a$.
\end{enumerate}

\begin{lemma}
\label{Newton Polygon}
Let $\ell= p^k$ with $k \geq 1$ and $d \in \overline{K}$. Let NP be the Newton polygon  of 
$$F(z)=(z+y)^\ell-y^\ell-d.$$ 
Let $m_1$ and $m_\ell$ be the first and last slope of NP respectively.
Set $\lambda_0 =\infty$,\ $\lambda_{k+1}=-\infty$, and  
$\lambda_n=k-n+\frac{p}{p-1}+\ell v(y)$ for $n = 1,\dots,k$.
If $d\neq 0$, choose $0\leq n_0 \leq k$ so that $\lambda_{n_0+1} \leq  v(d) < \lambda_{n_0}$. 

\begin{enumerate}
     \item The $x$-coordinates of the vertices of NP are exactly $ 0,  p^{n_0}, p^{n_0+1},\dots,p^k.$
\item If $n_0 \leq k-1$, then 
\[\displaystyle m_1= \frac{k-n_0+(\ell-p^{n_0}) v(y)-v(d)}{p^{n_0}}
\quad \text{and} \quad\ m_\ell = -\frac{p}{\ell(p-1)}-v(y).
\]     
    \item If $n_0=k$, then $\displaystyle m_1= m_\ell = -\frac{v(d)}{\ell}.$
    \end{enumerate}

If $d=0$, then the vertices of NP have $x$-coordinates $p, p^2, \dots, p^k$. The first slope is $m_1 = \infty$ and the last slope is $m_\ell = -\frac{p}{\ell(p-1)}-v(y)$.
\end{lemma}
\begin{proof}    
Write 
\[
F(z) = (z + y)^\ell - y^\ell - d = \sum\limits_{n=1}^{\ell} \binom{\ell}{n} y^{\ell - n} z^n - d.
\]
By Lemma~\ref{val_binomial}, the Newton Polygon (NP) is the lower convex hull of the points
\begin{equation}
\label{eq:NPpts}
\{P_0= (0, v(d)), P_n=(n, k - v(n) + (\ell - n) v(y)) \colon 1 \leq n \leq \ell \}.
\end{equation}
To determine the vertices of NP for the family of functions $
F(z) = (z + y)^\ell - y^\ell - d,$
as $v(d)$ and $v(y)$ vary, we shift these points by the linear function $y = v(y) x - \ell v(y)$ to obtain a new set of points:
\begin{equation}
\label{eq:convhull}
\{Q_0= (0, v(d) - \ell v(y)), Q_n= (n, k - v(n)) \colon 1 \leq n \leq \ell \}.
\end{equation}
As discussed earlier, the $x$-coordinates of the vertices of NP are precisely the $x$-coordinates of the vertices of the convex hull of these points. Note that the coordinates of all points, except for the first one in \eqref{eq:convhull}, are independent of $v(y)$ and $v(d)$.(See Example 
\ref{ex:NPeg} below for an explicit illustration of this construction.)
We first consider the points in ~\eqref{eq:convhull} excluding the $Q_0$: \begin{equation} 
\label{eq:convhull2}
\{Q_n =  (n, k - v(n)) \colon 1 \leq n \leq \ell \}.
\end{equation}
We claim that the lower convex hull of the points in \eqref{eq:convhull2} consists of  
$
Q_1, Q_p, \dots, Q_{p^j}, \dots, Q_{p^k}.
$
To prove this, it suffices to show that for any $p^{n-1} < i < p^n$, the point $Q_i$ lies above the line segment $Q_{p^{n-1}}Q_{p^n}$. The slope of $Q_{p^{n-1}}Q_{p^n}$ is  
$
\frac{-1}{p^n - p^{n-1}},
$  
and the slope of $Q_{p^{n-1}}Q_i$ is  
$
\frac{n-1 - v(i)}{i - p^{n-1}}.
$  
For $Q_i$ to lie above $Q_{p^{n-1}}Q_{p^n}$, we need  
$$
\frac{n-1 - v(i)}{i - p^{n-1}} > \frac{-1}{p^n - p^{n-1}}.
$$
This inequality holds because $p^{n-1} < i < p^n$ ensures $0 < i - p^{n-1} < p^n - p^{n-1}$, and $v(i) \leq n-1$ implies $n-1 - v(i) > -1$.
Thus, the convex hull of the points in \eqref{eq:convhull} consists of
$Q_0, Q_{p^j}, \dots, Q_{p^k}, \textrm{ for some } j = 0,1,\dots k.$
Here, $j$ depends on the $y$-coordinate of $Q_0$.
Switching back to the points in \eqref{eq:NPpts}, we see that the convex hull consists of: $P_0, P_{p^j}, \dots, P_{p^k}.$
Consider two points in the set of points \eqref{eq:NPpts} where the $x$-coordinates are consecutive powers of $p$: 
$
P_{p^{n-1}} = \left(p^{n-1},k-(n-1)+(\ell-p^{n-1})v(y)\right)
 \text{ and }  
P_{p^{n}}= \left(p^{n},k-n+(\ell-p^{n})v(y)\right).
$
The $y$-intercept of the line between these two points is given by: 
\[
\lambda_n=k-n+\frac{p}{p-1}+\ell v(y).
\]
Choose $n_0$ so that  $\lambda_{n_0} \geq  v(d) > \lambda_{n_0+1}$.  (Recall that we define $\lambda_{k+1} = -\infty$.)
\begin{itemize}
    \item If $n_0 <k$, we see that for all $i< n_0$, $P_{p^i}$ is above the line connecting $(0,v(d))$ and $P_{p^{n_0}}$, but $P_{p^{n_0}}$ is below the line connecting $(0,v(d))$ and $P_{p^{n_0+1}}$.
    \item  If $n_0=k$, then all points $P_{p^i}$ for $i<k$  lie above the line connecting $(0,v(d))$ and $P_{p^k}$, so NP is a single segment.
\end{itemize}
In both cases, we see that $j=n_0$ as desired. (Example \ref{ex: lamda} below illustrates the case $n_0 = 1$.)
The slopes of NP given in the lemma statement are a straightforward calculation from equation~\eqref{eq:NPpts} once we know the $x$-coordinates of NP.
In the case $d=0$, we use the same reasoning to see that NP has a vertical line to the point $\left(p, k-1+(\ell-p)v(b)\right)$. The slope of $m_1=\infty$ corresponds to the root at $z=0$. The vertices of the lower convex hull are then the points $(p^n, k-n + (\ell-p^n)v(b))$.
\end{proof}

\begin{example}
\label{ex:NPeg}
The following figure illustrates the shifted Newton Polygon of $(z+y)^8-y^8-d$. \\
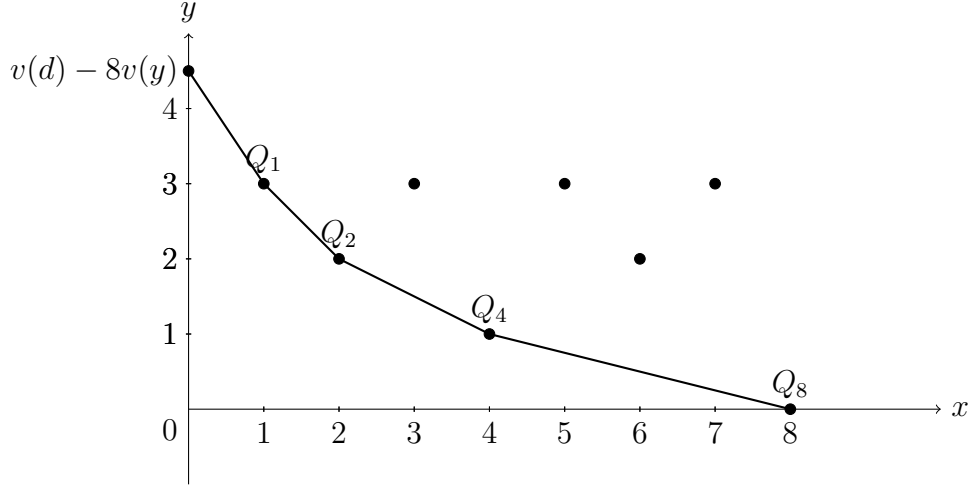
\begin{tikzpicture}
   
    \draw[->] (0,0) -- (10,0) node[right] {$x$};
   
    \draw[->] (0,-1) -- (0,5) node[above] {$y$};
   
    \draw[black, thick]  (0,4.5) -- (1,3) -- (2,2) -- (4,1) -- (8,0);
   
    \filldraw[black] (0,4.5) circle (2pt) node[above] {};    \filldraw[black] (1,1pt) -- (1,-1pt)  node[anchor=north] {$1$};
    \filldraw[black] (2,1pt) -- (2,-1pt)  node[anchor=north] {$2$};
    \filldraw[black] (3,1pt) -- (3,-1pt)  node[anchor=north] {$3$};
    \filldraw[black] (4,1pt) -- (4,-1pt)  node[anchor=north] {$4$};
    \filldraw[black] (5,1pt) -- (5,-1pt)  node[anchor=north] {$5$};
    \filldraw[black] (6,1pt) -- (6,-1pt)  node[anchor=north] {$6$};
    \filldraw[black] (7,1pt) -- (7,-1pt)  node[anchor=north] {$7$};
    \filldraw[black] (8,1pt) -- (8,-1pt)  node[anchor=north] {$8$}; \filldraw[black]  (0,1) node[left] {$1$};
    \filldraw[black] (0,2)   node[left] {$2$};
    \filldraw[black] (0,3)  node[left] {$3$};
    \filldraw[black]  (1pt,4) -- (-1pt,4);
    \filldraw[black] (0,4)  node[left] {$4$};
    \filldraw[black]  (1pt,1) -- (-1pt,1);
    \filldraw[black]  (1pt,2) -- (-1pt,2);
    \filldraw[black]  (1pt,3) -- (-1pt,3);
    \filldraw[black] (0,0) node[anchor=north east] {$0$};
    \filldraw[black] (0,1)  node[left] {$1$};
    \filldraw[black] (0,2)  node[left] {$2$};
    \filldraw[black] (0,3)  node[left] {$3$};
    \filldraw[black] (2,2) circle (2pt) node[above]{$Q_{2}$};
    \filldraw[black] (4,1) circle (2pt) node[above] {$Q_{4}$};
    \filldraw[black] (8,0) circle (2pt) node[above]{$Q_{8}$};
    \filldraw[black] (1,3) circle (2pt) node[above]{$Q_{1}$};
    \filldraw[black] (3,3) circle (2pt);
    \filldraw[black] (5,3) circle (2pt);
    \filldraw[black] (6,2) circle (2pt);
    \filldraw[black] (7,3) circle (2pt);

    \draw (0,4.5) node [left] {$v(d)-8v(y)$};
           
    \end{tikzpicture}
\end{example}

\begin{example}
\label{ex: lamda}
The following figure illustrates the Newton Polygon of $(z+y)^8-y^8-d$ with the corresponding $\lambda_n$, where $\lambda_1 > v(d) > \lambda_2$. 

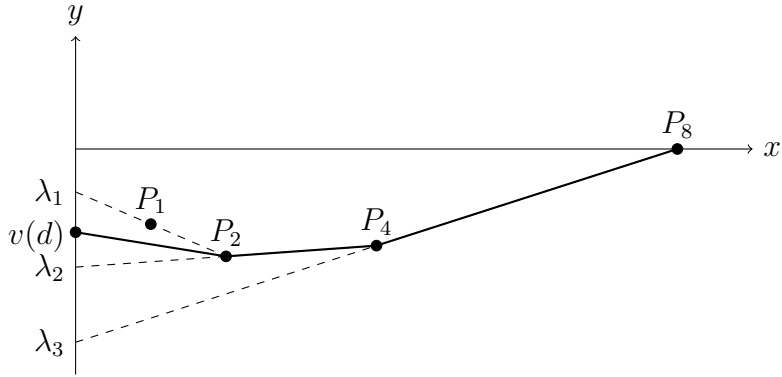
\begin{tikzpicture}
   
    \draw[->] (0,0) -- (9,0) node[right] {$x$};
   
    \draw[->] (0,-3) -- (0,1.5) node[above] {$y$};
   
    \draw[black, thick]  (0,-6/7-0.25) -- (2,-10/7) -- (4,-9/7) -- (8,0);
   
    \filldraw[black] (0,-6/7-0.25) circle (2pt) node[above] {};    \filldraw[black] (1,-1) circle (2pt) node[above] {$P_{1}$};
    \filldraw[black] (2,-10/7) circle (2pt) node[above]{$P_{2}$};
    \filldraw[black] (4,-9/7) circle (2pt) node[above] {$P_{4}$};
    \filldraw[black] (8,0) circle (2pt) node[above]{$P_{8}$};
   
    \draw (0,-6/7-0.3) node [left] {$v(d)$};
   
    \draw [dashed] (0,-4/7) node [left] {$\lambda_1$}-- (2,-10/7);
    \draw [dashed] (0,-11/7) node [left] {$\lambda_2$}-- (4,-9/7);    
    \draw [dashed] (0,-18/7) node [left] {$\lambda_3$}-- (8,0);        
    \end{tikzpicture}
\end{example}

For the remaining case of $\ell$, we have a similar result. The proof is essentially the same (just the last endpoint changes), so we omit the details.
\begin{lemma}
\label{lem: General NP}
Let $\ell=Np^k$ where $p\nmid N$ and $N>1$. Let NP be the Newton polygon  of $$F(z)=(z+y)^\ell-y^\ell-d.$$ Let $m_1$ and $m_\ell$ be the first and last slope of NP respectively.\\
Set 
$\lambda_n(y)=k-n+\frac{p}{p-1}+\ell v(y)$ for $n = 1,\dots,k$ and $\lambda_{k+1}(y) = \ell v(y).$ Then
\begin{enumerate}
\label{Newton Polygon 2}
    \item If $d=0$, then the vertices of NP have $x$-coordinates $p, p^2, \dots, p^k, Np^k$. The first slope is $m_1 = \infty$ and the last slope is $m_\ell = -v(y)$.
    \item If  $\lambda_{n_0+1}(y) <  v(d) \leq \lambda_{n_0}(y)$ for some integer $0\leq n_0 \leq k$, then the $x$-coordinates of the vertices of NP are exactly $0,  p^{n_0}, p^{n_0+1},\dots,p^k,Np^k,$
        and $$\displaystyle m_1= \frac{k-n_0+(\ell-p^{n_0}) v(y)-v(d)}{p^{n_0}} \ \textrm{and }\ m_l=  -v(y).  $$
    \item If $v(d) \leq \lambda_{k+1}(y)$, then $m_1=m_\ell=-\frac{v(d)}{\ell}$.
   \end{enumerate}
\end{lemma}

Similar to Lemma 3.1 in \cite {anderson-hamblen-poonen-walton}, we obtain the following corollary.
\begin{corollary}
\label{cor: root in base field}
     Let $d,y \in \overline{K}$ and let $f_c(z)=z^\ell-c$. Then
     \begin{enumerate}
         \item If $v(d)> k-1+\frac{p}{p-1}+\ell v(y)$, the equation $f_c(x)-f_c(y)=d$ has at least one root $t \in K(d,y)$.
         \item If $\ell = p^k$ where $k \geq 1$ and $v(d)< \frac{p}{p-1}+\ell v(y)$, all the roots $t$ of the equation $f_c(x)-f_c(y)=d$ have valuation $\frac{v(d)}{\ell}$.
         \item If $\ell=Np^k$ where $(N,p)=1$ and $N>1$, and  $v(d)< \ell v(y)$, all the roots $t$ of the equation $f_c(x)-f_c(y)=d$ have valuation $\frac{v(d)}{\ell}$.
     \end{enumerate}
\end{corollary}
\begin{proof} Let $z = x-y$ then $f_c(x)-f_c(y)-d = (z+y)^\ell-y^\ell-d$. The result follows from Lemma \ref{Newton Polygon} and Lemma \ref{Newton Polygon 2}.
\end{proof}

\section{Julia sets of $z^\ell - c$}\label{sec:JuliaSets}

In this section, we describe the Berkovich filled Julia set of the polynomial $f_c(z)=z^\ell -c$ for $\ell \geq 2$ following the ideas in \cite{Ben23}. See Chapter~6 for background on $\mathbb{P}^1_{\an}$, the Berkovich projective line.

First, we recall the definition of the {\it Weierstrass degree} of a power series on a closed disk, which is helpful in determining the image of a power series (in our case, the polynomial $f_c(z)$) on a closed disk. 
\begin{definition}
    Let $F(z) = \sum_{n \geq 0} c_n (z-a)^n$ be a nonzero power series converging on a closed disk $\bar D = \bar D (a,r)$ for $a \in \mathbb{C}_v$ and $r>0$. The Weierstrass degree $\wdeg_{\bar D}(F)$ is the largest integer $d \geq 0$ such that  $|c_d|r^d= \max_{n\geq 0}|c_n|r^n$.   
\end{definition}

\begin{theorem}[\protect{\cite[Theorem 3.15]{Ben23}}]
\label{thm: radius}
Let $\bar D \subset \mathbb{C}_v$ be a closed disk of radius $r>0$ containing a point $a \in \mathbb{C}_v$, and let $F(z)= \sum_{n \geq 0} c_n(z-a)^n$ be a nonconstant power series converging on $\bar D$. Let $d:=\wdeg_{\bar D}(F-c_0)$, and let $t:= |c_d|r^d>0.$ Then $F (\bar D)$ is a closed disk, centered at $c_0$, and of radius $t$. Moreover, $F:\bar D \longrightarrow F(\bar D)$ is everywhere $d$-to-$1$, counting multiplicity.     
\end{theorem}

\begin{lemma}[\protect{\cite[Lemma 14.8]{Ben23}}] 
\label{lem:preimdiskrad}
Let $\bar D, \bar E \subset\mathbb{C}_v$ be closed disks, let $a\in \bar E$, let $F$ be a power series converging on $\bar E$, and suppose that $F(a) \in \bar D \subset F(\bar E)$. Then there is a unique closed disk $\bar U \in \bar E$ with $a\in \bar U$ such that $F(\bar U)=\bar D$.
\end{lemma} 
\begin{corollary}
\label{cor: pre-radius}
    Let $F(z)= \sum_{i \geq 1}^n c_i z^i$ be a polynomial, and let $\bar D = \bar D(0,R)$ be a disk centered at $0$ with radius $R >0$. Then in $F^{-1}(\bar D)$, the disk $U$ containing $0$ has radius  $r= \left(\frac {R} {|c_{m}|}\right)^{1/m}$ where $m = \wdeg_{\bar U}(F)$.
\end{corollary}
\begin{proof} 
Let  $\bar E $ be a closed disk with $0 \in \bar E$ and radius $r_E$. Since $F(0)=0$, Theorem \ref{thm: radius} says $F(\bar E)$ is a closed disk containing 0 of radius $|c_m|r_E^m$ for some $m\geq 1$. So we may choose $r_E$ large enough that $\bar D \subseteq F(\bar E)$. By Lemma~\ref{lem:preimdiskrad}, there is a unique closed disk $\bar U \subseteq \bar E$ with $0 \in \bar U$ and $F(\bar U) = \bar D$. From Theorem~\ref{thm: radius}, the radius $r$ of $\bar U$ must satisfy $R = |c_m|r^m$.
\end{proof}

Let $b\in K$ be a fixed point of $f_c(z)$.
As in the proof of Theorem \ref{prop:goodreduction}, conjugate $f_c$ by a map that moves $b$ to 0  to obtain 
\begin{equation}
\label{series of g}
g(z)=(z+b)^\ell-b-c = \sum\limits_{n=1}^{\ell}a_n z^{n} \textrm{ where } a_n=\binom{\ell}{n}b^{\ell-n}.
\end{equation}
Since $g(0)= 0$, we have $b^\ell-b=c$, which will be useful in this section.

\begin{notation}
\label{radii}    
Let $r_0=0$; for $1\leq n\leq k$, let $r_n=p^{(-1/(p^{n}-p^{n-1}))-v(b)}$; and let \\
\[
r_{k+1}= \begin{cases}
     \infty &\text{ if } N=1 \text { and } k \geq 1, \\
    p^{-v(b)} & \text{ if } N>1.
\end{cases}
\] 

\end{notation}
We have the following lemma on the Weierstrass degree of $g(z)$ on disks of various radii. It follows from the notation above that Case~\eqref{N>1} never happens when $N=1$.

\begin{lemma} 
\label{lem: wdeg} 
Continue with the notation in \ref{radii} and let $r\in \mathbb{R}_{\geq 0}$. 
\begin{enumerate}
    \item\label{N=1} If there is some $n$ so that $r_{n}\leq r <r_{n+1},$ then $  \wdeg_{\bar D(0,r)}(g)=p^n$.  
    \item\label{N>1} If $r_{k+1} \leq r$, then $\wdeg_{\bar D(0,r)}=\ell$.
\end{enumerate}

\end{lemma}

\begin{proof}
The Newton polygon of $g(z)$ is described in Lemma~\ref{Newton Polygon}, using $y=b$ and $d=0$. A calculation shows that 
  the slope of the $n^\text{th}$ segment (between the points with $x$-coordinates $p^{n-1}$ and $p^n$) is given by $\log_{p}r_n$.

Hence for $r_{n}\leq r <r_{n+1}$, the disk $\bar D(0,r)$  contains exactly $p^n$ of the roots of $g(z)$. By Theorem \ref{thm: radius}, we conclude that $\wdeg_{\bar D(0,r)}(g)=p^n$. If $r_{k+1}\leq r$, all the roots of $g$ would be contained by $\bar D(0,r)$, hence $\wdeg_{\bar D(0,r)}(g)=\ell$.
\end{proof}

\begin{definition}
    Let $F$ be a polynomial of degree $d\geq 2$. We define the Berkovich filled Julia set of $F$ as
    $$ \mathcal{K}_{F,\an} := \{ \zeta  \in \mathbb{P}^1_{\an}: \lim_{n \to \infty} F^n(\zeta) \neq \infty  \}. $$
The Berkovich Julia set, denoted $\mathcal{J}_{F,\an}$, is  the boundary of $\mathcal{K}_{F,\an}.$ 
\end{definition}

\begin{prop}
    The Berkovich Julia set $\mathcal{J}_{g,\an}$ is the singleton set $\{\zeta(0,1)\}$ if and only if $v(c)\geq \nu_{\good}.$
\end{prop}
\begin{proof}
      By Proposition \ref{prop:goodreduction}, the hypothesis gives that $f_c$ has potential good reduction, so $g$ has good reduction.
By \cite[Proposition 8.12]{Ben23}, the Julia set for a map with good reduction is just the Gauss point.
\end{proof}

It remains to understand the Berkovich Julia set of $g$ when $g$ has bad reduction; that is when $v(c) < \nu_{\good}$. Let $R_0$ be the radius of the smallest disk centered at 0 containing all the roots of $g$. From the Newton polygon of $g$, we conclude: 
\[
R_0= 
\begin{cases}
     p^{(-1/(p^{k}-p^{k-1}))-v(b)}
    &\text{ if } N=1 \text{ and } k \geq 1, \\
    p^{-v(b)}& \text{ if }  N>1.
\end{cases}
\] 
By Lemma \ref{lem: wdeg},  for $|z|>R_0 $, $z^\ell$ is the dominant term in $g(z)$, so $|g(z)|=|z|^\ell.$ 
Since $v(c)<\nu_{\good} \leq 0$ and $b^\ell-b = c$, necessarily $v(b)<0$. It follows that $R_0 >1,$ so $\displaystyle \lim_{n \to \infty}g^{n}(z)=\infty$ for $|z|>R_0 $. 
Therefore, the Berkovich filled Julia set $\mathcal{K}_{g,\an}$ is contained in the closed disk $\bar{D}(0,R_0):= U_0 $. 

We consider the preimages of $U_0$, denoted as $U_m:= g^{-m}(U_0),$ for $m\geq 0.$
We prove the following properties of $U_m$. 
\begin{lemma}
\label{lem: radius of preimage}
Suppose $c\in \overline K$ with $v(c)<\nu_{\good}$.
\begin{enumerate}
    \item\label{part1} For each $\alpha_m \in g^{-m}(0)$, we have $v(\alpha_m+b) = v(b)$ and $v(\alpha_m) \geq v(b)$.
    \item\label{part2} Each $U_m$ is the disjoint union of finitely many disks of the same radius.
     \item\label{part3} $U_m$ is nested, i.e \[U_0 \supsetneq U_1\supsetneq U_2\supsetneq \dots.\]
    \item\label{part4} The intersection $\displaystyle \bigcap_{i \geq 0} U_i$ is the disjoint union of disks of radius $R$ for some $R \geq 0$.
\end{enumerate}

\end{lemma}
\begin{proof}
   We prove~\eqref{part1} and~\eqref{part2} by induction on $m$. Since $b^\ell-b=c$ and $v(c)<0$, we see that $v(b)=v(c)/ \ell > v(c)$.  

       For $\alpha_1 \in g^{-1}(0)$, since $g(\alpha_1)=(\alpha_1+b)^\ell-b-c =0$, then 
       \[v(\alpha_1+b) 
       \ =\ v(b+c)/\ell 
       \ \geq\ \min \{v(b),v(c)\}/\ell 
       \ =\ v(c)/\ell
       \ =\ v(b).\] 
       Thus $v(\alpha_1+b) =v(b)$ and $v(\alpha_1) \geq v(b)$.

    Now let $\alpha_m \in g^{-1}(\alpha_{m-1})$. By the induction hypothesis, 
    \[
    v(b) \leq v(\alpha_{m-1}) =v((\alpha_m+b)^\ell -b-c).
    \]
   Since $v(c)<v(b)$, 
 it follows that $v(\alpha_m+b) = v(c)/\ell$, so $v(\alpha_m+b)= v(b)$ and $v(\alpha_m) \geq v(b)$.

For statement~\eqref{part2}, since $g^{-m}(0)$ is a finite set, there will be finitely many disjoint disks in $g^{-m}(U_0)$. We need to show that the radii of these disks are all equal. For $m=0$, the statement is obvious. Now assume it is true for $m-1$. 

For all $m \geq 0$, let $R_m$ be the radius of the disk in $g^{-m}(U_0)$ centered at $0$, and let $\alpha_m \in g^{-m}(0)$.  We want to show that 

\[
g(\bar D(\alpha_m,R_m))=\bar D(g(\alpha_m), R_{m-1})).
\]
We can recenter the polynomial $g(z)$ around $\alpha_m$ 
\[g(z)= (z-\alpha_m + \alpha_m +b)^\ell-b-c= \displaystyle \sum_{n=1}^{\ell}\binom{\ell}{n} (\alpha_m +b)^{\ell-n}(z-\alpha_m)^n + g(\alpha_m).\] 
By part~\eqref{part1}, $v(\alpha_m+b)= v(b)$. Hence, $v\left (\binom{\ell}{n} (\alpha_m+b)^{\ell-n}\right)= v(a_n)$, where $a_n$ is defined as in equation~\eqref{series of g}. Let $c_0:=g(\alpha_m)$ and $c_n:=\binom{\ell}{n} (\alpha_m+b)^{\ell-n}$ for all $n\geq 1$.
We have \[|a_n|{R_m}^n=|c_n|{R_m}^n\] for all $n \geq 1$. Therefore, we have $
\wdeg_{\bar D(0,R_m)}(g(z)-a_0)=\wdeg_{\bar D(\alpha_m,R_m)}(g(z)-c_0)$, and by Theorem 4.2,
\[
g(\bar D(\alpha_m,R_m))=\bar D(g(\alpha_m), R_{m-1})).
\]
Thus in the preimage  $g^{-m}(U_0)$, the disks centered at 0 and $\alpha_m$ have the same radius. 

To prove statement~\eqref{part3}, it suffices to consider those disks centered at~0.
 Recall that $U_0=\bar D(0,R_0)$ with 
\[
v(R_0) = \begin{cases}
     -\frac{p}{\ell(p-1)} -\frac{v(c)}{\ell}
    &\text{ if } N=1 \text{ and } k \geq 1, \\
    -\frac{v(c)}{\ell}& \text{ if }  N>1.
\end{cases}
\]  
Let $D(0,r)$ be a disk in $U_1=g^{-1}(U_0))$ such that $g(D(0,r))= U_0.$
Assume that $r \geq R_0$, then by Lemma \ref{lem: wdeg}, $\textrm{wdeg}_{\bar D(0,r)}(g)=\ell$, hence $R_0= r^{\ell}$ by Theorem \ref{thm: radius}. 
The condition $v(c)<\nu_{\good}$ implies that $R_0>1$, then $R_0= r^{\ell} > r$. This contradiction shows that $U_1 \subsetneq U_0$.   
For each $k \geq 1$, we observe that $U_{k+1}=f_c^{-k}(U_1)  , U_k=f_c^{-k}(U_0),$ and $U_1 \subsetneq U_0$. Therefore, $U_{k+1} \subsetneq U_{k}.$ This proves~\eqref{part3}.

For part~\eqref{part4}, let $R_m$ denote the radius of each disk in $U_m$. By part~\eqref{part3}, the sequence $\{R_m\}$ is decreasing and convergent, with $R = \lim_{m \to \infty} R_m$. 
In the Berkovich space, every point of $\bigcap_{m \geq 1} U_m$ is either a type II point with radius $R$ or a type IV point. In particular, since $0 \in \bigcap_{m \geq 1} U_m$, we know $\bigcap_{m \geq 1} U_m$ contains $\overline{D}(0, R)$.
\end{proof}
The value of the radius $R$ depends on the value of $v(c)$. We use the following notation.
\begin{notation}
\label{level of c}
    Let  $c_0= -\infty$; for $n = 1,\dots, k$, let  $c_n:= \frac{-\ell}{\ell-1} (k-n+ \frac{p^{n}-1}{p^n-p^{n-1}})$; in the case $N>1$, let $c_{k+1}=0$.   

\end{notation}

Note that $c_0 < c_1 < \cdots < c_k < c_{k+1}$, and we have \[
\nu_{\good} = \begin{cases}
     c_k
    &\text{ if } N=1 , \\
    c_{k+1}& \text{ if }  N>1.
\end{cases}
\]   The values of $c_n$ are chosen due to the following lemma.
\begin{lemma}  
\label{lem: c level} 
 Let $\{a_n\}$ be the coefficients of $g(z)$ in equation~\eqref{series of g} and let $R$ be the radius in Lemma~\ref{lem: radius of preimage}\eqref{part4}.
Suppose $c_{n}\leq v(c)< c_{n+1}$ for 
some $n= 0,1,\dots, k$. 
\begin{enumerate}
    \item\label{case1} If $n \geq 1$, then  $R= |a_{p^{n}}|^{-1/(p^{n}-1)}$.
   
    \item\label{case2} If $n = 0$, then $R=0$. 
    
\end{enumerate}

\end{lemma}
\begin{proof} 

Define the function $T:\mathbb{R^{+}} \to \mathbb{R}\cup\{\infty\}$ by \[T(s) = ||g(z)||_{\bar D(0,s)} = \max_{n\geq 0}|a_n|s^n.\] 
Continuing with the notation in the proof of Lemma \ref{lem: radius of preimage}, let $R_m$ be the radius of each disjoint disk in $g^{-m}(U_0)$. By Corollary \ref{cor: pre-radius}, we have $$T(R_m)= R_{m-1}.$$ Since $T(s)$ is a continuous function, taking a limit on both sides yields $T(R)=R$.

Using Notation~\ref{radii}, let $n_0$ be the integer such that $r_{n_0}\leq R < r_{n_0+1}.$ Then 
$$R =T(R)=|a_{p^{n_0}}|R^{p^{n_0}} .$$
If $n_0 \neq 0$, then $R = |a_{p^{n_0}}|^{-1/(p^{n_0}-1)}$.

We prove that for $n\geq 1$ 
\begin{equation}
\label{level of limit}
    r_n \leq  |a_{p^{n}}|^{-1/(p^{n}-1)} < r_{n+1} \quad \textrm {if and only if} 
    \quad c_n \leq v(c) <c_{n+1},
\end{equation} 
which will prove part~\eqref{case1}.

The left inequalities in~\eqref{level of limit} are equivalent to 
\begin{equation}
    \label{eq:imptineq} v(r_n) \leq \frac{v(a_{p^n})}{p^n-1}   < v(r_{n+1}).
    \end{equation}
 
By Lemma \ref{val_binomial}, we have $$v(a_{p^n})= v\left(\binom{\ell}{p^n} b^{\ell-p^n}\right)= k-n +\frac{(\ell-p^n)}{\ell}v(c).$$
Combining this with the definition of $r_n$ and $r_{n+1}$ for $n<k$, equation~\eqref{eq:imptineq} becomes
\[
-\frac{p^n-1}{p^n-p^{n-1}}-\frac{p^n-1}{\ell}v(c) \leq 
 k-n +\frac{(\ell-p^n)}{\ell}v(c) 
< - \frac{p^n-1}{p^{n+1}-p^{n}} -\frac{p^n-1}{\ell}v(c).
\]

\[
-\frac{\ell}{\ell-1}\left(k-n+\frac{p^n-1}{p^n-p^{n-1}}\right)
\leq  
v(c) 
<
 - \frac{\ell}{\ell-1} \left( k-n+\frac{p^n-1}{p^{n+1}-p^{n}} \right).
\] 
By Notation \ref{level of c}, if $n<k$ this is equivalent to $c_n \leq v(c) <c_{n+1}$.

If $n=k$, the inequalities in~\eqref{eq:imptineq} become
\[-\frac{p^k-1}{p^k-p^{k-1}}-\frac{p^k-1}{\ell}v(c) \leq  \frac{(\ell-p^k)}{\ell}v(c) < -\frac{p^k-1}{\ell}v(c),\]
 which is equivalent to $c_k \leq v(c) < c_{k+1}=0.$

We are left with the case $n_0=0$, so that $R = |a_1|R$. We need to show that $|a_1| \neq 1$.
Note that if $v(c)\geq c_1$, then there is some $n_0>0$ such that $c_{n_0} \leq v(c) < c_{n_0+1}$, which is impossible by the chain of equivalences above.  So we must have $v(c) < c_1= -\frac{\ell k}{\ell-1}$. Therefore 
\[
v(a_1)= k-0 +\frac{(\ell-1)}\ell  v(c)
< k - k = 0.
\]
We conclude that  $|a_1|\neq 1$, and thus $R=0$.
\end{proof}

By Lemmas~\ref{lem: radius of preimage} and~\ref{lem: c level}, we deduce that the set  $\mathcal{K}_{g,\an}$ is homeomorphic to a Cantor set of type I points if and only if $n_0  = 0$. 
   
    In this case, we have  $\mathcal{K}_{g,\an}=\mathcal{J}_{g,\an}$. We sum up our results in the following theorem.

\begin{theorem}
\label{Thm: Julia sets}
Suppose that $\ell \in \mathbb{Z}$ and $c \in \bar K$. Set \[\nu_\infty = \displaystyle \frac{-\ell}{\ell-1}v(\ell) \quad \textrm{ and } \quad \nu_{\good} = 
\begin{cases}
     -\frac{p}{p-1}
    &\text{ if } \ell = p^k, \text{ for }k>1, \\
    0 & \text{ if } \ell = Np^k \text{ for }N>1
\end{cases}.\]
 \begin{enumerate}
     \item  If $v(c)<\nu_\infty$, then $\mathcal{J}_{g,\an}$ is homeomorphic to a Cantor set of type I points. 
     \item If $\nu_\infty\leq v(c)<\nu_{\good}$, then $\mathcal{J}_{g,\an}$  is homeomorphic to a Cantor set containing type II points and possibly type IV points.
     \item If $v(c) \geq \nu_{\good} $, then $\mathcal{J}_{g,\an}$ is a single type II point.  
 \end{enumerate}
\end{theorem}
\begin{remark}
    \begin{enumerate}
        \item When  $v(c) < \nu_\infty$ the Berkovich Julia set coincides with the classical Julia set. We show in the next section that it is a necessary and sufficient condition that the extension $K_\infty/K$ is finite. 
        \item When $v(c) \geq \nu_{\good}$, the polynomial $g(z)$ has good reduction as mentioned in a previous section.  
        \item In \cite{anderson-hamblen-poonen-walton}, the authors consider only the cases when $p \nmid \ell$ or $\ell=p$ (so $k = 1$). In these cases, the cutoff valuations coincide, i.e $\nu_\infty =\nu_{\good}$. 
    \end{enumerate}
    \end{remark}

\section {The extension $K_\infty$}\label{sec:extension}
In this section, we describe how the extension $K_{\infty}/K$ varies with $v(c)$. In~\cite[Corollary 4.4]{anderson-hamblen-poonen-walton}, the authors prove that if $p\nmid \ell$ or $p=\ell$ and $v(c) <\nu_\infty$, then $K_\infty$ is a finite extension of $K$.
The authors prove this by examining directly the valuation $v(\alpha_n)$ such that $f_c^n(\alpha_n)= \alpha$. 

Here we use the results of Section~\ref{sec:JuliaSets} to recover this result and extend to all values of $\ell \geq 2$. For simplicity, we first prove the case $\ell = p^k$. The general case $\ell = Np^k$ with $N>1$ and $(N,p)=1$ follows by the same argument, using Lemma~\ref{lem: General NP} in place of Lemma~\ref{Newton Polygon}. The proof proceeds identically except that the final segment of the Newton polygon is handled by the same method. The details are straightforward and left to the reader.

\subsection{Sufficiently Negative Values: $v(c) < v_{\infty}$}
We begin by extending the result that for sufficiently negative values of $v(c)$, the extension $K_\infty/K$ is in fact finite. For this, we prove that the Julia set $\mathcal{J}_{g,\an}$ is contained in a finite extension of $K$, and then we apply Krasner's Lemma \cite{NSW08} to deduce the same is true for $K_\infty$.

\begin{lemma}[Krasner's Lemma]
Let $K$ be a complete non-archimedean field and $K^{\text{sep}}$ a separable closure of $K$. Given an element $\alpha \in K^{\text{sep}}$, denote its Galois conjugates by $\alpha=\alpha_1, \alpha_2, \dots, \alpha_n$. Then if $\beta \in K^{\text{sep}}$ satisfies
$v(\alpha - \beta) > v(\alpha - \alpha_i)$ for $i\neq 1$ then $K(\alpha) \subseteq K(\beta)$.
\end{lemma}

\begin{prop}\label{thm:finiteJulia}
Suppose that $\ell \geq 2$ and $c \in \overline{K}$. 
    If $v(c)< \nu_\infty$, then $\mathcal{J}_{f,\an}$ is contained in a finite extension of $K$.
\end{prop}
\begin{proof}
Assume $v(c) < \nu_\infty$ and let $b$ be a fixed point of $f_c$. Since $g(z) = f_c(z+b) - b$, the Julia sets of $f$ and $g$ are just translates of each other. We consider the field extensions $K_n = K(f_c^{-n}(b))$, and $K_\infty = \cup_{n \geq 0} K_n$. (So here we take $\alpha = b$ to construct the field extension $K_\infty$.)By Theorem~\ref{Thm: Julia sets}(1), the Julia set $\mathcal{J}_{f_c,\an}$ contains only type I points. Let $R_m$ be the radius of the disk in $g^{-m}(U_0)$ centered at $0$, as in the proof of Lemma~\ref{lem: radius of preimage}. Since $\lim_{m \to \infty} R_m = 0$, we choose $n$ such that $R_n < r_1$. By Lemma~\ref{lem: wdeg}(1) and Lemma~\ref{lem: radius of preimage}, we have $\wdeg_{\bar{D}(\alpha_m, R_m)}(g) = 1$ for all $m \geq n$, $\alpha_m \in g^{-m}(0)$.  

We prove by strong induction that $K_m = K_n$ for all $m > n$. First consider the case $m=n+1$.  
Let $\alpha_{n+1} \in K_{n+1}$. Then $\alpha_{n+1} \in U_n$, and hence $\alpha_{n+1} \in \bar{D}({\alpha_n}', R_n)$ for some ${\alpha_n}' \in K_n$. Since $\wdeg_{\bar{D}({\alpha_n}', R_n)} = 1$, Theorem~\ref{thm: radius} implies $g$ maps $\bar{D}({\alpha_n}', R_n)$ one-to-one onto its image. Therefore, $\alpha_{n+1}$ is the unique root of the polynomial $g(z) - g(\alpha_{n+1}) \in K_n[z]$ in the disk $\bar{D}({\alpha_n}', R_n)$. By the uniqueness of $\alpha_{n+1}$, we have  
\[
v(\alpha_{n+1} - {\alpha_n}') > v(\alpha_{n+1} - \alpha_{n+1,i})
\]  
for all Galois conjugates $\alpha_{n+1,i}$ of $\alpha_{n+1}$ over the complete non-Archimedean field $K_n$. By Krasner's Lemma, we have $\alpha_{n+1} \in K_n({\alpha_n}') = K_n$.  

Let $m > n$ and assume that $K_m = K_{m-1} = \dots = K_n$. Let $\alpha_{m+1} \in K_{m+1}$. By the same argument above , we conclude that $\alpha_{m+1} \in $ $K_n$, and hence $K_{m+1} = K_n$.  
By induction, $K_m = K_n$ for all $m > n$, and thus $K_\infty = K_n$.  

Finally, by \cite[Proposition~5.23]{Ben23},  
\[
\mathcal{J}_{g,\an} = \overline{\cup_{n \geq 0} g^{-n}(0)} \subset \overline{K_n}.
\]  
Since $K$ is complete and $K_n$ is a finite extension of $K$, $K_n$ is complete, and hence $\overline{K_n} = K_n$. Therefore,  
\[
\mathcal{J}_{f_c,\an} \subset K_n(b),
\]  
which is a finite extension of $K$.  
       
\end{proof}

Proposition~\ref{thm:finiteJulia} says that if our root point $\alpha$ is in the Julia set of $f_c$, then $K_\infty/K$ is a finite extension. It remains to extend this to an arbitrary root point $\alpha \in K$.

\begin{theorem}\label{thm:suffnegc-pt1}
Suppose that $\ell \geq 2$ and $c \in \overline{K}$. 
    If $v(c)< \nu_\infty$, then $K_{\infty}/K$ is a finite extension.
\end{theorem}

\begin{proof}

Let $\alpha\in K$ and let $K_n=K(f_c^{-n}(\alpha))$. To show that $K_\infty = \cup_{n\geq0} K_n$ is a finite extension over $K$, it is sufficient for us to show that all the solutions of $g^{n}(z)=\alpha-b$ lie in $D(0,R_0)$ for sufficiently large $n$.

By Lemma 3.1, if $v(\alpha-b)<\lambda_{k}(b)$, then for any $x$ such that $g(x)=\alpha-b$, $v(x)=\frac{v(\alpha-b)}{\ell}$, hence $v(x)\geq\lambda_{k}(b)$ for sufficiently large $n$ and hence all solutions lie within $D(0,R_0)$, i.e.~all solutions of $f^{n}(z)=\alpha$ lie within $D(b,R_0)$.
\end{proof}
\subsection{Insufficiently Negative Values: $v_{\infty}\leq v(c)<0$} We recall the following notation from Section~\ref{sec:NP}:
\[
\lambda_n(y)=k-n+\frac{p}{p-1}+\ell v(y),
\]
which we will apply here in the case that $y=b$ (a fixed point of $f_c$)
so $v(y) = v(b) = v(c)/\ell$.

The following lemma says that  we can further assume that $v(\alpha_n)=v(c)/\ell$ for all $n\geq 1$ and all $\alpha_n \in f_c^{-n}(\alpha).$ 

\begin{lemma}[\protect{\cite[Lemma 3.2]{anderson-hamblen-poonen-walton}}]
\label{lem: v(alpha}
    Suppose that $v(c)<0$ and $\alpha_n \in f_c^{-n}(\alpha)$. If $n$ is sufficiently large, then  $v(\alpha_n) =v(c)/\ell.$ If $v(\alpha)>v(c)$, then this conclusion holds for all $n \geq 1.$
\end{lemma}

\begin{theorem}\label{thm:suffnegc}
Suppose that $\ell \geq 2$ and $\ell \neq p$. 
    \label{inf-ram}
    Then if $\nu_\infty\leq v(c)<0$, then $K_\infty / K$ is infinitely wildly ramified.

\end{theorem}

\begin{proof} 

    First, we begin with the case that $\nu_\infty\leq v(c)< -\frac{p}{p-1} $.  By Lemma $\ref{lem: c level}$ , we have $$\lim_{m\to\infty}R_m = |a_{p^{n}}|^{-1/(p^{n}-1)} \textrm{ for some }n\geq 1.$$

    Let $\alpha_0=\alpha$ and choose $\alpha_1 \in f^{-1}(\alpha_0)$ to be the furthest from $\alpha_0.$ Then for $n\geq 2$, we choose inductively $\alpha_{n+1} \in f_c^{-1}(\alpha_n)$ to be the closest to $\alpha_n$. Let $d_n= \alpha_{n+1}-\alpha_n.$ Since $v(c)<-\frac{p}{p-1}$, it follows that $\lambda_k(b) < 0$. We observe that if $v(d_0) \geq \lambda_k(b)$, then we are in Case~(2) of Lemma \ref{Newton Polygon}. We conclude that the furthest preimage $\alpha_1$ has valuation given by the first slope $m_\ell$, so 
    $$v(d_1) = \frac{p}{\ell(p-1)}+\frac{v(c)}{\ell}.$$
    Otherwise, if $v(d_0) < \lambda_k(b)$, then by Corollary \ref{cor: root in base field}, $v(d_1)=\frac{v(d_0)}{\ell}.$ Since $v(d_0) < \lambda_k(b)<0$, we see that $v(d_n) > \lambda_k(b)$ for $n$ large enough. Then $$v(d_{n+1})=\frac{p}{\ell(p-1)}+\frac{v(c)}{\ell}.$$
    By relabeling, we can assume that  $v(d_0) = \frac{p}{\ell(p-1)}+\frac{v(c)}{\ell}.$ Note that $v(d_0)= v(R_0)$ where $R_0$ is defined in Lemma \ref{lem: radius of preimage}.

    We will prove by induction that $v(d_n) = v(R_n) $ for all $n \in \mathbb{Z}.$ We already know $$v(d_0)=v(R_0)=\frac{p}{\ell(p-1)}+\frac{v(c)}{\ell}.$$ Assume that $v(d_m)=v(R_m)$ for all $0\leq m<n$. we make the following 2 claims: 

    \begin{enumerate}
        \item $-v(d_{n})$ is the slope of the first segment of the Newton Polygon of the polynomial $(z+\alpha_n)^\ell-{\alpha_n}^\ell-d_{n-1}$.
        \item $-v(R_n)$ is the slope of the first segment of the Newton Polygon of the polynomial $(z+b)^\ell-b^\ell-R_{n-1}$.
    
    \end{enumerate}
To prove (1), notice that for all $\beta \in f^{-1}(\alpha_n)$, the difference $\beta - \alpha_n$ satisfies the polynomial  
$$(z + \alpha_n)^\ell - \alpha_n^\ell - d_{n-1}.$$  
Since $\alpha_{n+1}$ is the root of $f_c(z) = \alpha_n$ closest to $\alpha_n$, the slope $-v(d_n)$ corresponds to the first segment of the Newton Polygon of the polynomial $(z + \alpha_n)^\ell - \alpha_n^\ell - d_{n-1}$.  
For (2), we write  
$$(z + b)^\ell - b^\ell - R_{n-1} = \sum\limits_{n=1}^{\ell} a_n z^n - R_{n-1}, \textrm{  where } a_n = \binom{\ell}{n}b^{\ell-n}.$$
The Newton Polygon of this polynomial is the lower convex hull of the points 
$$
\{(n, v(a_n)) : 1 \leq n \leq \ell\} \cup \{(0, v(R_{n-1}))\}.
$$
Therefore the slope of the first segment of the Newton Polygon is given by $$\min_{1\leq m\leq\ell}\frac{v(a_m)-v(R_{n-1})}{m}=-\max_{1\leq m\leq\ell}v((\frac{R_{n-1}}{a_m})^{\frac{1}{m}})=-v(R_n).$$
The last equality follows from Corollary \ref{cor: pre-radius}. 

    Since $v(\alpha_n)=v(b)=\frac{v(c)}{\ell}$ and $v(d_{n-1})=v(R_{n-1})$, the polynomials $(z+\alpha_n)^\ell-{\alpha_n}^\ell-d_{n-1}$ and $(z+b)^\ell-b^\ell-R_{n-1}$ have the same Newton Polygon and hence $v(d_n)=v(R_n).$
    
    Then by Lemma \ref{lem: c level},  we have
    \[\lim_{n\rightarrow\infty}v(d_n)=-\frac{v(a_{p^{n_0}})}{p^{n_0}-1} \textrm { for some } n_0 = 1,2, \dots, k. \] and there is an integer $m$ such that for $n\geq m$ we have \[\lambda_{n_0+1} \leq v(d_n) < \lambda_{n_0}.\]
    Setting $q=p^{n_0}$ and using Lemma \ref{Newton Polygon} , we compute $v(d_n)$ recursively and deduce that
    \begin{align*}
      v(d_{m+n}) = & \ \frac{v(d_m)}{q^n}-\left(\frac{1}{q}+\frac{1}{q^2}+\dots+\frac{1}{q^n}\right)v(a_q)\\
      =& \ \frac{v(d_m)}{q^n}-\frac{q^n-1}{q^n(q-1)}v(a_q).
    \end{align*}\\
    Since $\displaystyle \lim_{n \to \infty} v(d_{m+n}) = - \frac{a_q}{q-1},$ the exponent of $p$ in the denominator of $v(d_{m+n})$ in the reduced form is unbounded as $n$ increases, hence $K_\infty/K$ is infinitely wildly ramified.
    
For the case $-\frac{p}{p-1}\leq v(c) < 0.$ Consider two subcases:
    \begin{enumerate}
          \item If $v(d_0)\leq \lambda_k$, then by Corollary \ref{cor: root in base field}, it is easy to show that for any $n \geq 0$, we have  \[v(d_n) \leq \lambda_k \textrm { and }v(d_n)= \frac{v(d_0)}{\ell^n}.\]
        \item If $v(d_0) >\lambda_k$, then $v(d_1)= \frac{p}{\ell(p-1)}+\frac{v(c)}{\ell} \leq \frac{p}{p-1}+v(c) = \lambda_k$. Then in this case we have
        \[v(d_n)= \frac{v(d_1)}{\ell^{n-1}}, \textrm {  for } n \geq 1.\]  
    \end{enumerate}
      
    In both cases, the extension $K_\infty/K$ is infinitely wildly ramified.
\end{proof}

\begin{theorem}\label{thm:suffnegc-part2}
Suppose that $\ell \geq 2, \ (\ell,p) \neq 1,$ and $c \in \overline{K}$. 
\label{finite-ram}If $v(c)= \nu_\infty$, then $K_\infty/K$ is an infinite extension, and  it is finitely ramified if and only if $\ell=p$ and $\alpha$ lies in the closed unit disk centered at a fixed point of $f$.

\end{theorem}

\begin{proof}
The result for $k=1$ is proved in {\cite[Theorem 1.3 ]{anderson-hamblen-poonen-walton}}.
     Suppose $k\geq 2$, then $\lambda_k = \frac{p}{p-1}+v(c) =\frac{p}{p-1}-\frac{k\ell}{\ell-1}  < 0$. Using the same argument in Theorem \ref{thm:suffnegc}, we deduce that $K_\infty /K$ is infinitely wildly ramified for all $\alpha$.
\end{proof}

\subsection{Nonnegative Valuation: $v(c) \geq 0$}  
In this case, we follow the argument presented in \cite{anderson-hamblen-poonen-walton}. For completeness, we provide the proof below.  

\begin{theorem}\label{thm:nonnegc}  
Let $\ell \geq 2$ and $\ell \neq p$. If $v(c) \geq 0$, then the extension $K_\infty / K$ is infinitely wildly ramified.  
\end{theorem}  

\begin{proof}  
Fix a sequence $(\alpha_n)$ such that $\alpha_0 = \alpha$ and $\alpha_{n+1} \in f^{-1}(\alpha_n)$. Select $\beta_0 = \alpha$, and $\beta_{n+1} \in f^{-1}(\beta_n)$. We choose $\beta_1 \in f^{-1}(\beta_0)$ to be the element furthest from $\alpha_1$. Set $d_n = \beta_n - \alpha_n$.  

First we consider the case $v(\alpha) \neq v(c)$.  If $\min(v(\alpha), v(c)) \neq 0$, then $v(\alpha_1) = \min\{v(\alpha), v(c)\} /{\ell}$. By induction, $v(\alpha_n) = \min\{v(\alpha), v(c)\}/ \ell^n$. Since $(\ell, p) \neq 1$, it follows that $K_\infty / K$ is infinitely wildly ramified.  

    Now, suppose $\min(v(\alpha), v(c)) = 0$.  We first consider the cases $v(c) > v(\alpha) = 0$ and $v(\alpha) > v(c) = 0$. 

    If $v(c) > v(\alpha) = 0$, then  $v(\alpha_n) = v(\beta_n) = 0$ for all $n$. Applying Lemma~\ref{Newton Polygon} with $d = 0$ and $y = \beta_1$, we deduce that  
    \[
    v(d_1) = \frac{p}{\ell(p-1)} < \lambda_k.  
    \]  
    By induction, $v(d_n) < \lambda_k$ for all $n \geq 1$, and thus  
    \[
    v(d_n) = \frac{v(d_1)}{\ell^n} = \frac{p}{\ell^{n+1}(p-1)}.  
    \]  
    Consequently, $K_\infty / K$ is infinitely wildly ramified.  

    In the case $v(\alpha) > v(c) = 0$, we have $v(\alpha_1) = 0$. Relabeling, we may assume $v(\alpha) = 0$ and consider this case with the remaining case $v(\alpha) = v(c) \geq 0$.

    Assume $v(\alpha_n) = v(c)$ for all $n$. If this does not hold, we reduce to one of the previous cases. Since $d_0 = 0$, Lemma~\ref{Newton Polygon} implies  
    \[
    v(d_1) = \frac{p}{\ell(p-1)} + v(c) < \lambda_k = \frac{p}{p-1} + \ell v(c).  
    \]  
    By induction, $v(d_n) < \lambda_k$ for all $n \geq 1$. Hence,  
    \[
    v(d_n) = \frac{v(d_1)}{\ell^n} = \left(\frac{p}{\ell(p-1)} + v(c)\right)/\ell^n.  
    \]  
    It follows that $K_\infty / K$ is infinitely wildly ramified.  
\end{proof}  

\end{document}